\documentclass[11pt, a4paper]{amsart}
\usepackage{mathbbol}
\usepackage{amsmath, amsthm, amssymb}
\usepackage{txfonts, mathrsfs}
\usepackage[usenames]{color}
\usepackage{psfrag}
\usepackage{graphicx}

\newtheorem*{thm*}{Theorem}

\newtheorem{thm}{Theorem}[section]
\newcommand{\bt}{\begin{thm}}
\newcommand{\et}{\end{thm}}

\newtheorem{cor}[thm]{Corollary}
\newcommand{\bc}{\begin{cor}}
\newcommand{\ec}{\end{cor}}

\newtheorem{lem}[thm]{Lemma}
\newcommand{\bl}{\begin{lem}}
\newcommand{\el}{\end{lem}}

\newtheorem{prop}[thm]{Proposition}
\newcommand{\bp}{\begin{prop}}
\newcommand{\ep}{\end{prop}}

\newtheorem{defn}[thm]{Definition}
\newcommand{\bd}{\begin{defn}}      
\newcommand{\ed}{\end{defn}}

\newtheorem{rmrk}[thm]{Remark}
\newcommand{\br}{\begin{rmrk}}
\newcommand{\er}{\end{rmrk}}

\newtheorem{quest}[thm]{Question}
\newcommand{\bq}{\begin{quest}}
\newcommand{\eq}{\end{quest}}

\newtheorem{example}[thm]{Example}

\newcommand{\N}{\mathbb{N}}

\newcommand{\R}{\mathbb{R}}

\newdimen\vintkern\vintkern12pt
\def\vint{-\kern-\vintkern\int}

\newcommand{\hm}{{\mathcal H}}

\newcommand{\diam}{\operatorname{diam}}
\newcommand{\trace}{\operatorname{tr}}

\newcommand{\Area}{\operatorname{Area}}

\newcommand{\md}{\operatorname{md}}

\newcommand{\bdry}{\partial}

\newcommand{\Vol}{\operatorname{Vol}}

\newcommand{\jac}{{\mathbf J}}
\newcommand{\ap}{\operatorname{ap}}
\newcommand{\apmd}{\ap\md}

\begin{document}
\bibliographystyle{plain}

\title[]{Area minimizing discs in locally non-compact metric spaces}

\author[C.-Y. Guo]{Chang-Yu Guo}
\address
  {Department of Mathematics\\ University of Fribourg\\ Chemin du Mus\'ee 23\\ 1700 Fribourg, Switzerland}
\email{changyu.guo@unifr.ch}

\author[S. Wenger]{Stefan Wenger}

\address
  {Department of Mathematics\\ University of Fribourg\\ Chemin du Mus\'ee 23\\ 1700 Fribourg, Switzerland}
\email{stefan.wenger@unifr.ch}

\keywords{}

\date{\today}

\thanks{Research partially supported by Swiss National Science Foundation Grants 153599 and 165848.}

\begin{abstract}
We solve the classical problem of Plateau in every metric space which is $1$-complemented in an ultra-completion of itself. This includes all proper metric spaces as well as many locally non-compact metric spaces, in particular, all dual Banach spaces, some non-dual Banach spaces such as $L^1$, all Hadamard spaces, and many more. Our results generalize corresponding results of Lytchak and the second author from the setting of proper metric spaces to that of locally non-compact ones. We furthermore solve the Dirichlet problem in the same class of spaces. The main new ingredient in our proofs is a suitable generalization of the Rellich-Kondrachov compactness theorem, from which we deduce a result about ultra-limits of sequences of Sobolev maps.
\end{abstract}

\maketitle

\section{Introduction and statement of main results}

The classical problem of Plateau is concerned with the existence of surfaces of disc-type with minimal area and prescribed Jordan boundary in Euclidean space. This problem has a long and rich history for which we refer for example to \cite{Dierkes-et-al10}. The first rigorous solution for arbitrary Jordan curves in Euclidean space was given independently by Douglas \cite{Dou31} and Rad\'o \cite{Rad30}. This solution was extended to a large class of Riemannian manifolds by Morrey \cite{Mor48}. Recently, Lytchak and the second author solved the classical problem of Plateau in the setting of arbitrary proper metric spaces in \cite{LW15-Plateau}. Recall that a metric space is proper if all its closed bounded subsets are compact. Before the paper \cite{LW15-Plateau} only a few results beyond the setting considered by Morrey existed, see \cite{Nik79}, \cite{MZ10}, \cite{OvdM14}. The existence and regularity results proved in \cite{LW15-Plateau} have had applications to problems in several fields, see \cite{LW-param}, \cite{LW-isoperimetric}, \cite{LWY16}. The purpose of the present note is to solve the classical Plateau problem as well as the Dirichlet problem in a class of metric spaces which includes many that are not locally compact. 

For a bounded domain $\Omega\subset \R^n$ with $n\geq 2$, a complete metric space $X$, and $p>1$ we denote by $W^{1,p}(\Omega, X)$ the space of Sobolev maps from $\Omega$ to $X$ in the sense of Reshetnyak \cite{Res97}. The Reshetnyak energy of $u\in W^{1,p}(\Omega, X)$ is denoted $E_+^p(u)$. If $\Omega\subset\R^n$ is a bounded Lipschitz domain then $u\in W^{1,p}(\Omega, X)$ has a trace, written as s$\trace(u)$, which belongs to $L^p(\partial \Omega, X)$. Let $D$ be the open unit disc in $\R^2$. The parametrized Hausdorff area of $u\in W^{1,2}(D, X)$ is denoted by $\Area(u)$. In \cite{LW15-Plateau} the authors introduced a notion of $Q$-quasiconformality for maps $u\in W^{1,2}(D, X)$. This is similar to but different from the notion of quasiconformal maps used in geometric function theory. We refer to Section~\ref{sec:prelims} for the definitions related to Sobolev maps mentioned above and for references.

We turn to our main results and first introduce the class of spaces for which we can solve the classical Plateau problem. We refer to Section~\ref{sec:prelims} for the notion of a non-principal ultrafilter $\omega$ on $\N$ and the definition of the ultra-limit $\lim\nolimits_\omega a_m$ of a bounded sequence $(a_m)$ of real numbers. Let $(X,d)$ be a metric space and $\omega$ a non-principal ultrafilter on $\N$. Denote by $X_\omega$ the set of equivalence classes $[(x_m)]$ of sequences $(x_m)$ in $X$ satisfying $\sup_m d(x_1, x_m)<\infty$, where sequences $(x_m)$ and $(x'_m)$ are identified if $\lim\nolimits_\omega d(x_m, x'_m) = 0$. The metric space obtained by equipping $X_\omega$ with the distance $d_\omega([(x_m)], [(x'_m)]) = \lim\nolimits_\omega d(x_m, x'_m)$ is called the ultra-completion or ultra-product of $X$ with respect to $\omega$. Clearly, $X$ isometrically embeds into $X_\omega$ via the map $\iota\colon X\hookrightarrow X_\omega$ which assigns to $x$ the equivalence class $[(x)]$ of the constant sequence $(x)$. 

\bd\label{def:1-compl-ultralimit}
A metric space $X$ is said to be $1$-complemented in some ultra-completion of $X$ if there exists a non-principal ultrafilter $\omega$ on $\N$ for which there is a $1$-Lipschitz retraction from $X_\omega$ to $X$. 
\ed

Our first main result can now be stated as follows. Given a Jordan curve $\Gamma\subset X$ we denote by $\Lambda(\Gamma, X)$ the possibly empty family of maps $v\in W^{1,2}(D, X)$ whose trace $\trace(v)$ has a continuous representative which is a weakly monotone parametrization of $\Gamma$.

\bt\label{thm:Plateau-1-complemented}
 Let $X$ be a complete metric space and $\Gamma$ a Jordan curve in $X$ such that $\Lambda(\Gamma, X)\not= \emptyset$. If $X$ is $1$-complemented in some ultra-completion of $X$ then there exists $u\in \Lambda(\Gamma, X)$ such that $$\Area(u) = \inf\{\Area(v): v\in\Lambda(\Gamma, X)\}$$ and $u$ is $\sqrt{2}$-quasiconformal.
\et

The class of spaces which are $1$-complemented in some ultra-completion includes all proper metric spaces, all dual Banach spaces, some non-dual Banach spaces such as $L^1$, furthermore all Hadamard spaces and injective metric spaces, see Proposition~\ref{prop:1-complemented-spaces}. Recall that a Hadamard space is a complete metric space which is ${\rm CAT}(0)$, that is, has non-positive curvature in the sense of Alexandrov. Our theorem thus applies to all these spaces and, in particular, for example to all $L^p$ spaces. It generalizes \cite[Theorems 1.1 and 10.2]{LW15-Plateau} and \cite{Nik79}. Exactly as in \cite{LW15-Plateau}, the quasiconformality constant $\sqrt{2}$ is optimal but can be improved to $1$ for metric spaces satisfying a certain property (ET) which, roughly speaking, excludes non-Euclidean normed spaces as weak tangents. 

We record the following special case of Theorem~\ref{thm:Plateau-1-complemented}.

\bc\label{cor:1-complemented-area}
 Let $X$ be a Hadamard space or a dual Banach space and $\Gamma\subset X$ a rectifiable Jordan curve. Then there exists $u\in \Lambda(\Gamma, X)$ such that $$\Area(u) = \inf\{\Area(v): v\in\Lambda(\Gamma, X)\}$$ and $u$ is $\sqrt{2}$-quasiconformal.
\ec

Since Hadamard spaces and Banach spaces admit a quadratic isoperimetric inequality in the sense of \cite{LW15-Plateau}, the regularity results in \cite{LW15-Plateau} imply that any $u$ as in the corollary has a locally H\"older (resp.~Lipschitz in the case that $X$ is a Hadamard space) continuous representative which extends continuously to the boundary $S^1$.

We can solve the Dirichlet problem in the same class of metric spaces:

\bt\label{thm:Dirichlet}
Let $X$ be a complete metric space, $\Omega\subset\R^n$ a bounded Lipschitz domain, and $w\in W^{1,p}(\Omega, X)$ for some $p>1$. If $X$ is $1$-complemented in some ultra-completion of $X$ then there exists $u\in W^{1,p}(\Omega, X)$ with $\trace(u) = \trace(w)$ and such that $$E_+^p(u) = \inf\left\{E_+^p(v): \text{ $v\in W^{1,p}(\Omega, X)$ and $\trace(v) = \trace(w)$}\right\}.$$
\et

The theorem furthermore holds with the Reshetnyak energy $E_+^p(u)$ replaced by the Korevaar-Schoen Dirichlet energy $E^p(u)$ defined in \cite{KS93}. The theorem generalizes for example \cite[Theorem 2.3]{LW16-harmonic} and \cite[Theorem 2.2]{KS93}. For regularity results for solutions of Dirichlet's problem in the metric space setting we refer for example to \cite{KS93} and \cite{LW16-harmonic} and the references therein.


The main new ingredient in the proofs of our results above is the following generalization of the Rellich-Kondrachov compactness theorem from the setting of proper metric spaces to that of arbitrary metric spaces.  Given a complete metric space $(X,d)$, a sequence $(u_m)\subset W^{1,p}(\Omega, X)$ will be called bounded if $$\sup_{m\in\N}\left[ \int_\Omega d(x_0, u_m(z))^p\, dz + E_+^p(u_m)\right]<\infty$$ for some and thus every $x_0\in X$.  For $p>1$ and $n\in\N$ we define the Sobolev conjugate of $p$ by $p^* = \frac{np}{n-p}$ if $p<n$ and $p^*=\infty$ otherwise.

\bt\label{thm:gen-Rellich-intro}
 Let $X$ be a complete metric space, $\Omega\subset \R^n$ a bounded Lip\-schitz domain, and $(u_m)\subset W^{1,p}(\Omega, X)$ a bounded sequence for some $p>1$. Then, after possibly passing to a subsequence, there exist a complete metric space $Z$, isometric embeddings $\varphi_m\colon X \hookrightarrow Z$ and $v\in W^{1,p}(\Omega, Z)$ such that $\varphi_m\circ u_m$ converges to $ v$ in $L^q(\Omega, Z)$ for every $q<p^*$. 
\et

A more general statement, which applies to sequences of metric spaces, will be given in Section~\ref{sec:Rellich}. Convergence in $L^q(\Omega, Z)$ means that $$\int_\Omega d(\varphi_m\circ u_m(z), v(z))^q\,dz \to 0$$ as $m\to\infty$. The proof of the theorem is similar to that for proper metric spaces but, in addition, uses a variant of Gromov's compactness theorem for sequences of metric spaces. Since the limit map $v$ can be viewed as a map to the ultra-completion $X_\omega$ for any $\omega$, see Lemma~\ref{lem:map-into-ultracompletion}, we will obtain the following result on ultra-limits of (sub)sequences of Sobolev maps.

\bt\label{thm:ultra-limit-Sobolev-maps}
 Let $X$ be a complete metric space and $X_\omega$ an ultra-completion. Let $\Omega\subset\R^n$ be a bounded Lipschitz domain and $(u_m)\subset W^{1,p}(\Omega, X)$ a bounded sequence for some $p>1$. Then, after possibly passing to a subsequence, the map $\phi(z):= [(u_m(z))]$ belongs to $W^{1,p}(\Omega, X_\omega)$ and satisfies $$E_+^p(\phi)\leq \liminf_{m\to\infty} E_+^p(u_m).$$ Moreover, if $\trace(u_m)$ converges to some map $\rho\in L^p(\partial \Omega, X)$ almost everywhere on $\partial \Omega$ then $\trace(\phi) = \iota\circ \rho$. Finally, if $p\geq n$ then $$\Vol_{\mu}(\phi)\leq \liminf_{m\to\infty} \Vol_\mu(u_m)$$ for any definition of volume $\mu$ (in the sense of convex geometry) inducing quasi-convex $n$-volume densities. 
\et

We refer to Section~\ref{sec:prelims} and \cite{LW15-Plateau} for the definition of the parametrized $\mu$-volume $\Vol_\mu(u)$ of a map $u\in W^{1,n}(\Omega, X)$ and for the notion of definition of volume inducing quasi-convex $n$-volume densities. Here, we simply mention that the Hausdorff measure provides a definition of volume which induces quasi-convex $2$-volume densities by \cite{BI12}. Moreover, in this case, the parametrized $\mu$-volume of a map $u\in W^{1,2}(D, X)$ coincides with the the parametrized Hausdorff area $\Area(u)$ used in Theorem~\ref{thm:Plateau-1-complemented} and Corollary~\ref{cor:1-complemented-area}. Theorem~\ref{thm:Plateau-1-complemented} and its corollary actually hold with the parametrized Hausdorff area replaced by the parametrized $\mu$-volume induced by any definition of volume $\mu$ which induces quasi-convex $2$-volume densities. 

We finally mention that sometimes it is possible to solve Plateau's problem even if a space is not $1$-complemented in an ultra-completion of itself. Indeed,
the Banach space $c_0$ of sequences of real numbers converging to $0$, equipped with the sup-norm, is not $1$-complemented in any ultra-completion of $c_0$. Nevertheless, the Plateau problem in $c_0$ has a solution for every Jordan curve $\Gamma\subset c_0$ for which $\Lambda(\Gamma, c_0)\not=\emptyset$, in particular whenever $\Gamma$ is rectifiable, see Proposition~\ref{prop:c_0-Plateau-not-1-complemented}. We do not know whether there exists a Banach space $X$ and a rectifiable Jordan curve in $X$ which does not bound an area minimizer. In contrast, it is not difficult to construct a complete metric space for which the Plateau problem is not solvable for some rectifiable Jordan curve and which, in addition, admits a quadratic isoperimetric inequality in the sense of \cite{LW15-Plateau}. See Example~\ref{ex:complete-metric-no-Plateau} below.

\bigskip

The paper is structured as follows. In Section~\ref{sec:prelims} we recall the necessary definitions concerning ultra-completions of metric spaces and Sobolev maps from a Euclidean domain to a complete metric space. In Section~\ref{sec:Rellich} we prove Theorem~\ref{thm:gen-Rellich} which implies Theorem~\ref{thm:gen-Rellich-intro}. Section~\ref{sec:ultra-limit-Sobolev-maps} is devoted to the proof of Theorem~\ref{thm:ultra-limit-Sobolev-maps}. In Sections~\ref{sec:Plateau} and \ref{sec:Dirichlet} we combine Theorem~\ref{thm:ultra-limit-Sobolev-maps} with the arguments in \cite{LW15-Plateau} and \cite{LW16-harmonic} to prove Theorems~\ref{thm:Plateau-1-complemented} and \ref{thm:Dirichlet}. In Section~\ref{sec:Plateau} we furthermore obtain an analog of Theorem~\ref{thm:Plateau-1-complemented} with area replaced by energy. We also provide Proposition~\ref{prop:c_0-Plateau-not-1-complemented} and Example~\ref{ex:complete-metric-no-Plateau} to which we alluded above.

\bigskip

{\bf Acknowledgments:} We are indebted to the referee for useful comments and questions which led to the statement of Theorem~\ref{thm:ultra-limit-Sobolev-maps}.

\section{Preliminaries}\label{sec:prelims}

\subsection{Basic notation}

The Euclidean norm of a vector $v\in\R^n$ is denoted by $|v|$ and the open unit disc in $\R^2$  by $D = \{v\in \R^2: |v|<1\}$. We write $\overline{D}$ for the closure of $D$ and $S^1$ for its boundary. 

Let $X$ be a metric space. A Jordan curve in $X$ is a subset of $X$ homeomorphic to $S^1$. A curve of finite length is called rectifiable. Let $\Gamma\subset X$ be a Jordan curve. A map $\gamma\colon S^1\to \Gamma$ is called a weakly monotone parametrization of $\Gamma$ if $\gamma$ is a uniform limit of homeomorphisms $\gamma_i\colon S^1\to \Gamma$.
For $m\geq 0$ the $m$-dimensional Hausdorff measure on $X$ is denoted by $\hm^m$. The normalizing constant is chosen in such a way that $\hm^m$ coincides with the $m$-dimensional Lebesgue measure on Euclidean $\R^m$. The Lebesgue measure of a set $A\subset \R^m$ is denoted by $|A|$.

\subsection{Ultra-completions of metric spaces}

We briefly recall the relevant definitions concerning ultra-completions and ultra-limits of metric spaces. We refer for example to \cite{BrH99} for details.

A non-principal ultrafilter on $\N$ is a finitely additive probability measure $\omega$ on $\N$ such that every subset of $\N$ is measurable and such that $\omega(A)$ equals $0$ or $1$ for all $A\subset \N$ and $\omega(A)=0$ whenever $A$ is finite. 
Given a compact Hausdorff topological space $(Z, \tau)$ and a sequence $(z_m)\subset Z$ there exists a unique point $z_\infty\in Z$ such that $\omega(\{m\in\N : z_m\in U\}) = 1$
for every $U\in\tau$ containing $z_\infty$. We denote the point $z_\infty$ by $\lim\nolimits_\omega z_m$.

Let $(X,d)$ be a metric space and $\omega$ a non-principal ultrafilter on $\N$. A sequence $(x_m)\subset X$ is called bounded if $\sup_m d(x_1, x_m)<\infty$. Define an equivalence relation $\sim$ on bounded sequences in $X$ by considering $(x_m)$ and $(x'_m)$ equivalent if $\lim_\omega d(x_m, x'_m) = 0$. Denote by $[(x_m)]$ the equivalence class of $(x_m)$. The ultra-completion $X_\omega$ of $X$ with respect to $\omega$ is the metric space given by the set $$X_\omega:= \left\{[(x_m)]: \text{ $(x_m)$ bounded sequence in $X$}\right\},$$ equipped with the metric $$d_\omega([(x_m)], [(x'_m)]):= \lim\nolimits_\omega d(x_m, x'_m).$$ Ultra-completions are sometimes called ultra-products in the literature. The ultra-completion $X_\omega$ of $X$ is the ultra-limit of the constant sequence $(X, x_0)$ with respect to $\omega$ for some fixed $x_0\in X$. 

The ultra-completion $X_\omega$ of $X$ is a complete metric space, even if $X$ itself is not complete. Notice that $X$ isometrically embeds into $X_\omega$ via the map $\iota\colon X\hookrightarrow X_\omega$ given by $\iota(x):= [(x)]$.

We now show that the classes of metric spaces mentioned after Theorem~\ref{thm:Plateau-1-complemented} satisfy Definition~\ref{def:1-compl-ultralimit}. 

\bp\label{prop:1-complemented-spaces}
 The class of metric spaces $X$ which are $1$-complemented in every ultra-completion of $X$ includes:
 \begin{enumerate}
  \item Proper metric spaces.
  \item Hadamard spaces.
  \item Injective metric spaces.
  \item Dual Banach spaces.
  \item Banach spaces which are $1$-complemented in some dual Banach space.
 \end{enumerate}
\ep

A metric space $X$ is said to be $1$-complemented in some metric space $Y$ if $X$ isometrically embeds into $Y$ and if there exists a $1$-Lipschitz retraction from $Y$ to $X$. This explains the terminology used in (v). Particular examples of  spaces satisfying (v) are given by $L$-embedded Banach spaces, see \cite{HWW93}, which includes $L^1$-spaces. A metric space $X$ is called injective if $X$ is $1$-complemented in every metric space into which $X$ embeds isometrically. We refer for example to \cite{Lan13} for properties of injective metric spaces.

\begin{proof}
Let $X$ be a metric space and $\omega$ a non-principal ultrafilter on $\N$. If $X$ is proper then the map $\iota$ is surjective. In particular, $X$ is $1$-complemented in $X_\omega$, which proves (i). If $X$ is a Hadamard space then so is $X_\omega$ and the orthogonal projection from $X_\omega$ to $X$ is $1$-Lipschitz since $X$ is a closed convex subset of $X_\omega$, see \cite{BrH99}. This proves (ii). If $X$ is an injective metric space then $X$ is $1$-complemented in $X_\omega$ by the definition of injectivity. This yields (iii). Let now $X$ be a dual Banach space. Closed balls of finite radius in $X$ are weak$^*$-compact by the Banach-Alaoglu theorem and the norm on $X$ is weak$^*$-lower semi-continuous. Thus, the map $P\colon X_\omega\to X$  given by $P([(x_n)]):= \lim\nolimits_\omega x_n$ is well-defined and $1$-Lipschitz. This proves (iv). The same argument works for (v).
\end{proof}

We end this subsection with the following easy observation which will be used in the proof of Theorem~\ref{thm:ultra-limit-Sobolev-maps}. It shows that the limit map $v$ appearing in Theorem~\ref{thm:gen-Rellich-intro} can be viewed as a map to an ultra-completion of $X$.

\bl\label{lem:map-into-ultracompletion}
Let $A$ be a set, $X$ a metric space, and $X_\omega$ an ultra-completion of $X$. Let $x_0\in X$ and let $f_m\colon A \to X$ be maps, $m\in\N$. Suppose there exist a metric space $Z$, isometric embeddings $\varphi_m\colon X\hookrightarrow Z$, and a map $g\colon A \to Z$ such that  $(\varphi_m(x_0))$ is a bounded sequence in $Z$ and $\varphi_m\circ f_m$ converges to $g$  pointwise on $A$. Then the map $\psi\colon g(A) \to X_\omega$ given by $\psi(g(a)):= [(f_m(a))]$ for $a\in A$ is well-defined and an isometric embedding.
\el

\begin{proof}
 We denote the distance on $X$ and $Z$ by $d$ and $d_Z$, respectively. Fix $z_0\in Z$ and notice that $\sup_m d_Z(z_0, \varphi_m(x_0))<\infty$. 
 For $a\in A$ we have $$d(x_0, f_m(a)) = d_Z(\varphi_m(x_0), \varphi_m(f_m(a))) \leq d_Z(z_0, \varphi_m(x_0)) + d_Z(z_0, \varphi_m(f_m(a)))$$ for all $m$. Since the right-hand side in the above inequality is bounded in $m$ it follows that $(f_m(a))$ is a bounded sequence in $X$.
Let $a, a'\in A$. Then $$d(f_m(a), f_m(a')) = d_Z(\varphi_m(f_m(a)), \varphi_m(f_m(a'))) \to d_Z(g(a), g(a'))$$ as $m\to\infty$, which implies that $\psi$ is well-defined and an isometric embedding.
\end{proof}

\subsection{Metric space valued Sobolev maps}

We briefly review the main definitions concerning Sobolev maps from a Euclidean domain to a metric space which will be used in the present paper.  We refer for example to \cite{LW15-Plateau} for details. There exist several equivalent definitions of Sobolev maps from Euclidean domains with values in a metric space, see e.g.~\cite{Amb90}, \cite{KS93}, \cite{Res97}, \cite{HKST15}, \cite{AT04}. Here, we recall the definition from \cite{Res97} using compositions with real-valued Lip\-schitz functions.

Let $(X,d)$ be a complete metric space and $p>1$ and let $\Omega\subset \R^n$ be a bounded domain. We denote by $L^p(\Omega, X)$ the set of measurable and essentially separably valued maps $u\colon \Omega\to X$ such that for some and thus every $x\in X$ the function $u_x(z):= d(x, u(z))$ belongs to $L^p(\Omega)$, the classical space of $p$-integrable functions on $\Omega$. A sequence $(u_m)\subset L^p(\Omega, X)$ is said to converge in $L^p(\Omega, X)$ to a map $u\in L^p(\Omega, X)$ if $$\int_\Omega d(u_m(z), u(z))^p\,dz \to 0$$ as $m\to\infty$.
 A map $u\in L^p(\Omega, X)$ belongs to the Sobolev space $W^{1,p}(\Omega, X)$ if there exists $g\in L^p(\Omega)$ such that for every $x\in X$ the function $u_x$ belongs to the classical Sobolev space $W^{1,p}(\Omega)$ and has weak gradient bounded by $|\nabla u_x|\leq g$ almost everywhere. 
The Reshetnyak $p$-energy of $u\in W^{1,p}(\Omega, X)$ is defined by $$E_+^p(u):= \inf\left\{\|g\|_{L^p(\Omega)}^p\;\big|\; \text{$g$ as above}\right\}.$$ There exist other natural definitions of energy of a Sobolev map, for example the well-known Korevaar-Schoen Dirichlet energy $E^p(u)$ defined in \cite{KS93}.

Let $\Omega\subset\R^n$ be a bounded Lipschitz domain. The trace of a Sobolev map $u\in W^{1,p}(\Omega, X)$ is defined as follows. Let $J=(-1,1)$ and $I=(-1,0)$. Given $x\in\bdry\Omega$ there exists an open neighborhood $U\subset\R^n$ of $x$, an open set $V\subset\R^{n-1}$, and a biLipschitz homeomorphism $\varphi\colon V\times J\to U$ of $x$ such that $\varphi(V\times I) = U\cap\Omega$ and $\varphi(V\times\{0\}) = U\cap \bdry\Omega$. For $\hm^{n-1}$-almost every $v\in V$ the map $t\mapsto u\circ\varphi(v,t)$ is in $W^{1,p}(I, X)$ and thus has an absolutely continuous representative, again denoted by $u\circ\varphi(v,\cdot)$. For $\hm^{n-1}$-almost every point $z\in U\cap\bdry\Omega$ the trace of $u$ at $z$ is defined by $$\trace(u)(z):= \lim_{t\to0^{-}} u\circ\varphi(v, t),$$ where $v\in V$ is such that $\varphi(v,0) = z$. By \cite[Lemma 1.12.1]{KS93} the definition of $\trace(u)$ is independent of the choice of $\varphi$. Since $\partial \Omega$ can be covered by a finite number of biLipschitz maps it follows that $\trace(u)$ is well-defined $\hm^{n-1}$-almost everywhere on $\bdry\Omega$. Furthermore, $\trace(u)$ is in  $L^p(\bdry\Omega, X)$  by \cite[Theorem 1.12.2]{KS93}, the definition of $L^p(\bdry\Omega, X)$ being analogous to that of $L^p(\Omega, X)$.

As was shown in \cite{Kar07} and \cite{LW15-Plateau}, every Sobolev map $u\in W^{1,p}(\Omega, X)$ has an approximate metric derivative at almost every point $z\in \Omega$ in the following sense. There exists a unique seminorm on $\R^n$, denoted $\apmd u_z$, such that 
 \begin{equation*}
    \ap\lim_{z'\to z}\frac{d(u(z'), u(z)) - \apmd u_z(z'-z)}{|z'-z|} = 0,
 \end{equation*}
where $\ap\lim$ denotes the approximate limit, see \cite{EG92}. If $u$ is Lipschitz then the approximate limit can be replaced by an honest limit.

Recall from \cite{AlvT04}  that a definition of volume (in the sense of convex geometry) is a function $\mu$ that assigns to every $n$-dimensional normed space $V$ a norm $\mu_V$ on $\Lambda^n V$ in such a way that $\mu_V$ is induced by the Lebesgue measure if $V$ is Euclidean and such that for every linear $1$-Lipschitz map $T\colon V\to W$ between $n$-dimensional normed spaces $V$ and $W$ the induced map $T_*\colon \Lambda^nV\to\Lambda^nW$ is also $1$-Lipschitz. Well-known examples of definitions of volume are the Busemann definition, the Holmes-Thompson definition, and the Benson definition. The Busemann definition is exactly the definition of volume which gives rise to the Hausdorff $n$-measure on $V$. The Benson definition is sometimes called the Gromov mass$^*$ measure.

\bd
 The parametrized $\mu$-volume of a map $u\in W^{1,n}(\Omega, X)$ is defined by $$\Vol_\mu(u):= \int_\Omega \jac^\mu(\apmd u_z)\,dz,$$ where the $\mu$-Jacobian $\jac^\mu(s)$ of a seminorm $s$ on $\R^n$ is given by $$\jac^\mu(s):= \left\{\begin{array}{l@{ \quad }l}
 \mu_{(\R^n, s)}(e_1\wedge\dots\wedge e_n) & \text{$s$ is a norm}\\
 0 & \text{otherwise}
 \end{array}\right.$$ 
and $e_1,\dots,e_n$ denote the standard basis vectors in $\R^n$.
\ed

The $\mu$-volume and Reshetnyak energy of $u\in W^{1,n}(\Omega, X)$ are related by $$\Vol_\mu(u)\leq E_+^n(u)$$ for any definition of volume $\mu$, see \cite[Lemma 7.2]{LW15-Plateau}. A definition of volume $\mu$ is said to induce quasi-convex $n$-volume densities if for any finite dimensional normed space $Y$ and any linear map $L\colon \R^n\to Y$ we have $\Vol_\mu(L|_B)\leq \Vol_\mu(\psi|_B)$ for every smooth immersion $\psi\colon B\to Y$ with $\psi|_{\partial B} = L|_{\partial B}$, where $B$ is the closed unit ball in $\R^n$. The Busemann definition of volume induces quasi-convex $2$-volume densities by the recent result \cite{BI12}. The Benson definition induces quasi-convex $n$-volume densities for every $n\geq 1$, see \cite{AlvT04}.

If $n=2$ and $\mu$ is the Busemann definition of volume we will denote $\Vol_\mu(u)$ by $\Area(u)$ and call it the (parametrized) Hausdorff area of $u$. It follows from the area formula \cite{Kir94}, \cite{Kar07} that if $u\in W^{1,2}(D, X)$ satisfies Lusin's property (N) then $$\Area(u)=  \int_X\#\{z : u(z) = x\} \,d\hm^2(x).$$

We finally need the following notion of quasiconformality introduced in \cite{LW15-Plateau}.

\bd
 A map $u\in W^{1,2}(D, X)$ is called $Q$-quasiconformal if for almost every $z\in D$ we have $\apmd u_z(v)\leq Q\cdot\apmd u_z(w)$ for all $v,w\in S^1$.
\ed

If $u$ is $Q$-quasiconformal then $E_+^2(u) \leq Q^2\cdot \Area(u)$, see \cite[Lemma 7.2]{LW15-Plateau}.

\section{Generalized Rellich-Kondrachov theorem}\label{sec:Rellich}

In this section, we prove the following theorem which generalizes Theorem~\ref{thm:gen-Rellich-intro}. 

\bt\label{thm:gen-Rellich}
 Let $\Omega\subset\R^n$ be a bounded Lipschitz domain. For every $m\in\N$, let $(X_m, d_m)$ be a complete metric space, $K_m\subset X_m$ compact, and $u_m\in W^{1,p}(\Omega, X_m)$.  Suppose that $(K_m, d_m)$ is a uniformly compact sequence and
\begin{equation}\label{eq:bdd-energy-distance-seq-um}
 \sup_{m\in\N}\left[ \int_\Omega d_m(x_m, u_m(z))^p\, dz + E_+^p(u_m)\right]<\infty
 \end{equation}
 for some and thus every $x_m\in K_m$. Then, after possibly passing to a subsequence, there exist a complete metric space $Z$, isometric embeddings $\varphi_m\colon X_m \hookrightarrow Z$, a compact subset $K\subset Z$ and $v\in W^{1,p}(\Omega, Z)$ such that $\varphi_m(K_m)\subset K$ for all $m\in \N$ and $\varphi_m\circ u_m$ converges to $v$ in $L^q(\Omega, Z)$ for every $q<p^*$. 
\et

Recall that a sequence of compact metric spaces $(B_m, d_m)$ is called uniformly compact if $\sup_m\diam B_m<\infty$ and if for every $\varepsilon>0$ there exists $N\in\N$ such that every $B_m$ can be covered by at most $N$ balls of radius $\varepsilon$. 

The proof of the theorem is similar to that of \cite[Theorem 5.4.3]{AT04} but uses, in addition, the following variant of Gromov's compactness theorem for sequences of metric spaces established in \cite[Proposition 5.2]{Wen11-cpt}.

\bp\label{prop:extended-gromov-compactness-sets}
Let  $(X_m, d_m)$ be a sequence of metric spaces and, for each $m\in\N$, subsets
\begin{equation*}
 B_m^1\subset B_m^2\subset B_m^3\subset\dots\subset X_m.
\end{equation*}
If for every $k\in\N$ the sequence $(B_m^k, d_m)$ is uniformly compact then, after possibly passing to a subsequence, there exist a complete metric space $Z$, isometric embeddings $\varphi_m\colon X_m\hookrightarrow Z$ and compact subsets $Y^1\subset Y^2\subset\dots\subset Z$ such that $\varphi_m(B_m^k)\subset Y^k$ for all $m\in\N$ and $k\in\N$.
\ep

We turn to the proof of Theorem~\ref{thm:gen-Rellich} and fix $m\in\N$. Since $X_m$ embeds isometrically into an injective metric space we may assume that $X_m$ is itself injective. Indeed, every metric space $X$ isometrically embeds into the Banach space $\ell^\infty(X)$ of bounded functions on $X$ with the supremum norm and $\ell^\infty(X)$ is injective. Now, there exists a non-negative function $h_m\in L^p(\Omega)$ such that $\|h_m\|_{L^p(\Omega)}^p\leq C\cdot E_+^p(u_m)$ for some constant $C$ only depending on $\Omega$ and $p$ and such that
\begin{equation*}
 d_m(u_m(z), u_m(z')) \leq |z-z'|\cdot (h_m(z) + h_m(z'))
\end{equation*}
for all $z,z'\in \Omega$, see e.g.~\cite[Proposition 3.2]{LW15-Plateau} and its proof. For $k\in \N$ set $$A_m^k:=\{ z\in \Omega : h_m(z) \leq k\}$$ and notice that the restriction of $u_m$ to $A_m^{k}$ is $2k$-Lipschitz.

\bl\label{lem:bdd-diam}
There exist $k_0\in\N$ and $\lambda>0$ such that $u_m(A_m^k)\subset B(K_m, \lambda k)$ and $A_m^{k}\not=\emptyset$ for all $m\in\N$ and $k\geq k_0$.
\el

Here, $B(K_m, \lambda k)$ denotes the set of all $x\in X$ for which there exists $y\in K_m$ with $d(x,y)<\lambda k$.

\begin{proof}
 For each $m\in\N$, fix $x_m\in K_m$ and define $C_m^k:= \{z\in\Omega : d_m(x_m, u_m(z))\leq k\}$. By Chebyshev's inequality and \eqref{eq:bdd-energy-distance-seq-um} there exists $M>0$ such that
 \begin{equation}\label{eq:bdd-Amk}
  |\Omega\setminus A_m^k|\leq k^{-p}\int_\Omega h_m^p(z)\,dz\leq M\cdot k^{-p}
\end{equation} and $|\Omega\setminus C_m^k|\leq M\cdot k^{-p}$ for all $m$ and $k$. Thus, there exists $k_0\in\N$ such that $A_m^k\cap C_m^k\not=\emptyset$ for all $m\in\N$ and all $k\geq k_0$. Fix $z_0\in A_m^k\cap C_m^k$. Then for every $z\in A_m^k$ we have $$d_m(x_m, u_m(z)) \leq d_m(x_m, u_m(z_0)) + d_m(u_m(z_0), u_m(z)) \leq k+ 2k\diam(\Omega),$$ so the lemma follows. 
\end{proof}

Let $m\in\N$ and $k\geq k_0$. Since $X_m$ is injective there exists a $2k$-Lipschitz map $u_m^k\colon\overline{\Omega}\to X_m$ which agrees with $u_m$ on $A_m^k$.  We define for each $m\in\N$ an increasing sequence of subsets $B_m^{k_0}\subset B_m^{k_0+1}\subset \dots \subset X_m$ by $$B_m^k:= K_m\cup u_m^{k_0}(\overline{\Omega})\cup\dots\cup u_m^k(\overline{\Omega}).$$ For fixed $k\geq k_0$, the sequence of metric spaces $(B_m^k, d_m)$ is uniformly compact by Lemma~\ref{lem:bdd-diam} and since $u_m^j$ is $2j$-Lipschitz on the compact set $\overline{\Omega}$. 
Thus, by Proposition~\ref{prop:extended-gromov-compactness-sets} there exists, after possibly passing to a subsequence, a complete metric space $(Z, d_Z)$, isometric embeddings $\varphi_m\colon X_m\hookrightarrow Z$, and compact subsets $Y^{k_0}\subset Y^{k_0+1}\subset\dots\subset Z$ such that $\varphi_m(B_m^k)\subset Y^k$ for all $m$ and $k\geq k_0$. In particular, for every $m\in\N$ the set $\varphi_m(K_m)$ is contained in the compact set $K:=Y^{k_0}$. Moreover, the maps $v_m = \varphi_m\circ u_m$ belong to $W^{1,p}(\Omega, Z)$ and satisfy
\begin{equation}\label{eq:bounded-vj-inZ}
 \sup_{m\in\N}\left[ \int_\Omega d_Z(z_0, v_m(z))^p\, dz + E_+^p(v_m)\right]<\infty
\end{equation}
for some and thus every $z_0\in Z$.

\bl \label{lem:vmj-conv-L1}
There exists a subsequence $(v_{m_j})$ which converges in $L^1(\Omega, Z)$ to some $v\in L^1(\Omega, Z)$.
\el

\begin{proof}
For given $k\geq k_0$, the map $v_m^k:= \varphi_m\circ u_m^k$ is $2k$-Lipschitz and has image in the compact set $Y^k$ for every $m\in\N$. Thus, by the Arzel\`a-Ascoli theorem and by a diagonal sequence argument, there exist  integers $1\leq m_1<m_2<\dots$ such that, for every $k\geq k_0$, the sequence $(v_{m_j}^k)$ converges uniformly on $\Omega$ as $j\to\infty$. Lemma~\ref{lem:diff-um-umk} below shows that there exists $M>0$ such that
\begin{equation*}
\int_\Omega d_Z(v_{m_j}(z), v_{m_l}(z))\,dz \leq 2M\cdot k^{1-p} + \int_\Omega d_Z(v_{m_j}^k(z), v_{m_l}^k(z))\,dz
\end{equation*}
for all $j, l\in\N$ and every $k\geq k_0$. Hence, the integral on the left-hand side converges to $0$ as $j,l\to\infty$. This proves that $(v_{m_j})$ is a Cauchy sequence in $L^1(\Omega, Z)$ and hence that $v_{m_j}$ converges in $L^1(\Omega, Z)$ to some $v\in L^1(\Omega, Z)$.
\end{proof}

The following lemma was used in the proof above.

\bl\label{lem:diff-um-umk}
There exists $M>0$ such that 
 \begin{equation}
  \int_\Omega d_Z(v_m(z), v_m^k(z))\,dz \leq M\cdot k^{1-p}
 \end{equation}
 for all $m\in\N$ and every $k\geq k_0$.
\el

\begin{proof}
Let $z_0\in Z$. By Lemma~\ref{lem:bdd-diam}, the definition of $v_m^k$, and since $\varphi_m(K_m)\subset K$, there exists $M'>0$ such that $d_Z(z_0, v_m^k(z))\leq M'\cdot k$ for every $z\in \Omega$ and all $m\in\N$ and $k\geq k_0$. This together with H\"older's inequality and \eqref{eq:bdd-Amk} yields
\begin{equation*}
 \begin{split}
  \int_\Omega d_Z(v_m(z), v_m^k(z))\,dz &= \int_{\Omega\setminus A_m^k} d_Z(v_m(z), v_m^k(z))\,dz\\
  &\leq M'\cdot k\cdot  \left|\Omega\setminus A_m^k\right| + \int_{\Omega\setminus A_m^k} d_Z(z_0, v_m(z))\,dz\\
  &\leq M''\cdot k^{1-p} + \left|\Omega\setminus A_m^k\right|^{1-\frac{1}{p}}\cdot \left(\int_\Omega d_Z(z_0, v_m(z))^p\,dz\right)^{\frac{1}{p}}\\
  &\leq M'''\cdot k^{1-p}
 \end{split}
\end{equation*}
for some constants $M''$ and $M'''$ which do not depend on $m$ and $k$.
\end{proof}

\bl\label{lem:conv-in-Lq}
For every $q<p^*$ the maps $v_{m_j}$ and $v$ belong to $L^q(\Omega, Z)$ and the sequence $(v_{m_j})$ converges to $v$ in $L^q(\Omega, Z)$.
\el

\begin{proof} 
Fix $z_0\in Z$ and let $q<\bar{q}<p^*$. By \eqref{eq:bounded-vj-inZ} and the classical Sobolev embedding theorem, the real-valued functions $z\mapsto d_Z(z_0, v_{m_j}(z))$ belong to $L^{\bar{q}}(\Omega)$ and form a bounded sequence in $L^{\bar{q}}(\Omega)$. Since a subsequence of $(v_{m_j})$ converges to $v$ almost everywhere it follows with Fatou's lemma that $v\in L^{\bar{q}}(\Omega, Z)$ and hence $$L:= \sup_{j\in\N} \left[\int_\Omega d_Z(v_{m_j}(z), v(z))^{\bar{q}}\,dz\right]<\infty.$$   Let $\varepsilon>0$. Then the set $F_\varepsilon^j:= \{z\in \Omega: d_Z(v_{m_j}(z), v(z)) >\varepsilon\}$ satisfies $|F_\varepsilon^j|\to 0$ as $j\to\infty$ because, by Chebyshev's inequality, $$|F_\varepsilon^j| \leq \varepsilon^{-1} \cdot \int_\Omega d_Z(v_{m_j}(z), v(z))\,dz$$ for every $j\in\N$ and because $v_{m_j}$ converges to $v$ in $L^1(\Omega, Z)$ by Lemma~\ref{lem:vmj-conv-L1}. By H\"older's inequality, 
\begin{equation*}
 \int_{\Omega} d_Z(v_{m_j}(z), v(z))^q\,dz \leq  \varepsilon^q\cdot |\Omega| + \int_{F_\varepsilon^j} d_Z(v_{m_j}(z), v(z))^q\,dz
 \leq \varepsilon^q\cdot |\Omega| + L^{\frac{q}{\bar{q}}}\cdot |F_\varepsilon^j|^{1-\frac{q}{\bar{q}}}
\end{equation*}
and hence
\begin{equation*}
 \int_{\Omega} d_Z(v_{m_j}(z), v(z))^q\,dz\to 0
\end{equation*}
as $j\to \infty$. This shows that $v_{m_j}$ converges to $v$ in $L^q(\Omega, Z)$, completing the proof. 
\end{proof}

Lemma~\ref{lem:conv-in-Lq} implies that $v_{m_j}$ converges to $v$ in $L^p(\Omega, Z)$. Since $E_+^p(v_{m_j})$ is uniformly bounded in $j$ it thus follows from \cite[Theorem 1.6.1]{KS93} that $v\in W^{1,p}(\Omega, Z)$. This concludes the proof of Theorem~\ref{thm:gen-Rellich}.

Theorem~\ref{thm:gen-Rellich-intro} is a direct consequence of Theorem~\ref{thm:gen-Rellich}. The proof of the latter moreover shows the following:

\br\label{rem:extension-intro-Rellich}
 The isometric embeddings $\varphi_m\colon X\hookrightarrow Z$ in Theorem~\ref{thm:gen-Rellich-intro} can be chosen with the following additional property. Given compact sets $C_1\subset C_2\subset \dots \subset X$ there exist compact sets $Y^1\subset Y^2\subset \dots \subset Z$ such that $\varphi_m(C_k) \subset Y^k$ for all $m$ and $k$. 
\er

Indeed, in the proof of Theorem~\ref{thm:gen-Rellich} one simply defines the subsets $B_m^k$ by $$B_m^k:= C_k\cup u_m^{k_0}(\overline{\Omega})\cup\dots\cup u_m^k(\overline{\Omega})$$ for $k\geq k_0$ and sets $Y^k:= Y^{k_0}$ for $k\leq k_0$. The rest of the proof remains unchanged.

\section{Ultra-limits of subsequences of Sobolev maps}\label{sec:ultra-limit-Sobolev-maps}

In this section we prove Theorem~\ref{thm:ultra-limit-Sobolev-maps}. For this let $X$, $X_\omega$, $\Omega$, and $(u_m)\subset W^{1,p}(\Omega, X)$ be as in the statement of the theorem. Let $\mu$ be a definition of volume inducing quasi-convex $n$-volume densities. After possibly passing to a subsequence, we may assume that $$E_+^p(u_m)\to \liminf_{k\to\infty} E_+^p(u_k)$$ as $m\to\infty$ and, if $p\geq n$ then also $\Vol_\mu(u_m) \to \liminf_{k\to\infty}\Vol_\mu(u_k)$.

We apply Theorem~\ref{thm:gen-Rellich} and fix $x_0\in X$.  After possibly passing to a subsequence, there thus exist a complete metric space $(Z, d_Z)$, a compact subset $K\subset Z$, and isometric embeddings $\varphi_m\colon X\hookrightarrow Z$ such that $\varphi_m(x_0)\in K$ for all $m$ and such that $v_m:= \varphi_m\circ u_m$ converges in $L^p(\Omega, Z)$ to some $v\in W^{1,p}(\Omega, Z)$. After possibly passing to a further subsequence, we may assume that $v_m$ converges to $v$ almost everywhere on $\Omega$. 
 Let $N\subset \Omega$ be a subset of Lebesgue measure zero such that $v_m(z)$ converges to $v(z)$ for all $z\in \Omega\setminus N$. 
 
 Define a subset $B\subset Z$ by $$B:= \{v(z): z\in \Omega\setminus N\}.$$ Then the map $\psi\colon B\to X_\omega$ given by $\psi(v(z)) = [(u_m(z))]$ whenever $z\in \Omega\setminus N$ is well-defined and an isometric embedding by Lemma~\ref{lem:map-into-ultracompletion}. 
 Since $X_\omega$ is complete, the map $\psi$ extends to an isometric map on $\overline{B}$, which we denote by $\psi$ again. After redefining $v$ on a set of measure zero, we may assume that $v$ has image in $\overline{B}$ and so $v$ is an element of $W^{1,p}(\Omega, \overline{B})$. The map $$\phi(z):= \psi(v(z)) = [(u_m(z))]$$ then belongs to $W^{1,p}(\Omega, X_\omega)$ and, by the lower semi-continuity of the Reshetnyak energy \cite[Corollaries 5.7]{LW15-Plateau}, furthermore satisfies 
 \begin{equation}\label{eq:phi-energy-lower}
  E_+^p(\phi) = E_+^p(v)\leq \liminf_{m\to \infty} E_+^p(v_m) = \liminf_{m\to \infty} E_+^p(u_m).
 \end{equation} In case $p\geq n$, the lower semi-continuity of the $\mu$-volume \cite[Corollaries 5.8]{LW15-Plateau} implies that 
  \begin{equation}\label{eq:phi-vol-lower}
  \Vol_\mu(\phi) = \Vol_\mu(v) \leq \liminf_{m\to\infty}\Vol_\mu(v_m) = \liminf_{m\to\infty}\Vol_\mu(u_m).
 \end{equation} 
 This proves the first and third part of the theorem. 
 
We now suppose that, in addition, the sequence of traces $\trace(u_m)$ converges to some map $\rho\in L^p(\partial\Omega, X)$ almost everywhere on $\partial \Omega$. Since $\rho$ is measurable and essentially separably valued it follows from Lusin's theorem that there exist compact sets $A_1\subset A_2\subset\dots\subset\partial \Omega$ such that the restriction $\rho|_{A_k}$ is continuous for every $k\in\N$ and $\hm^{n-1}(\partial \Omega\setminus A_k)\to 0$. Thus, the sets $C_k:= \rho(A_k)$ are compact and satisfy $C_1\subset C_2\subset\dots\subset X$. 
Now, Theorem~\ref{thm:gen-Rellich} and Remark~\ref{rem:extension-intro-Rellich} show that, after possibly passing to a subsequence, there exist a complete metric space $(Z, d_Z)$, compact subsets $Y^1\subset Y^2\subset \dots\subset Z$, isometric embeddings $\varphi_m\colon X\hookrightarrow Z$ and $v\in W^{1,p}(\Omega, Z)$ such that $\varphi_m(C_k)\subset Y^k$ for all $m$ and $k$ and $v_m:= \varphi_m\circ u_m$ converges in $L^p(\Omega, Z)$ to $v$ as $m\to \infty$. 

Set $C:= \bigcup_{k=1}^\infty C_k$. After passing to a further subsequence we may assume that $v_m$ converges to $v$ almost everywhere on $\Omega$ and that $\varphi_m|_C$ converges pointwise to an isometric embedding $\varphi\colon C\hookrightarrow Z$, the convergence being uniform on each $C_k$. Let $N\subset \Omega$ be a set of Lebesgue measure zero such that $v_m(z)$ converges to $v(z)$ for all $z\in \Omega\setminus N$. 

Define a subset of $Z$ by $$B:= \{v(z): z\in\Omega\setminus N\}\cup \varphi(C).$$ The map $\psi\colon B\to X_\omega$ given by $\psi(v(z))= [(u_m(z))]$ when $z\in \Omega\setminus N$ and by $\psi(\varphi(x))= \iota(x)=  [(x)]$ when $x\in C$, is well-defined and an isometric embedding by Lemma~\ref{lem:map-into-ultracompletion}. Since $X_\omega$ is complete there exists a unique isometric extension of $\psi$ to $\overline{B}$, which we denote again by $\psi$. After possibly redefining the map $v$ on $N$, we may assume that $v$ has image in $\overline{B}$ and hence $v$ is an element of $W^{1,p}(\Omega, \overline{B})$. The map $$\phi(z):= \psi(v(z))= [(u_m(z))]$$ then belongs to $W^{1,p}(\Omega, X_\omega)$ and satisfies \eqref{eq:phi-energy-lower} and, if $p\geq n$, then also \eqref{eq:phi-vol-lower}. Moreover, we have that $\trace(v_m) = \varphi_m\circ\trace(u_m)$ converges to $\varphi\circ\rho$ almost everywhere on $\partial \Omega$ and a subsequence of $\trace(v_m)$ converges to $\trace(v)$ almost everywhere on $\partial \Omega$ by \cite[Theorem 1.12.2]{KS93}. It thus follows that $\trace(v) = \varphi\circ\rho$ and hence $$\trace(\phi) = \psi\circ\trace(v) = \psi\circ\varphi\circ\rho = \iota\circ\rho,$$ completing the proof of the second part of the theorem.

\section{Area and energy minimizers with prescribed Jordan boundary}\label{sec:Plateau}

In this section we prove Theorem~\ref{thm:Plateau-1-complemented} as well as an analog for the energy, see Theorem~\ref{thm:energy-min-1-complemented} below. The proofs of these theorems are almost the same as in the case of proper metric spaces, see \cite[Theorems 7.1 and 7.6]{LW15-Plateau}, but they make essential use of Theorem~\ref{thm:ultra-limit-Sobolev-maps} instead of the Rellich-Kondrachov compactness theorem for proper metric spaces.

Let $X$ and $\Gamma$ be as in the statement of Theorem~\ref{thm:Plateau-1-complemented}.

\bl\label{lem:keep-area-lower-energy}
 For every $v\in \Lambda(\Gamma, X)$ there exists $u\in \Lambda(\Gamma, X)$ which is $\sqrt{2}$-quasi\-conformal and satisfies $\Area(u)\leq \Area(v)$. 
\el

\begin{proof}
Let $v\in \Lambda(\Gamma, X)$ and define $\Lambda_v:= \{w\in\Lambda(\Gamma, X): \Area(w)\leq \Area(v)\}$, which is not empty. Let $(v_m)\subset\Lambda_v$ 
be an energy minimizing sequence in $\Lambda_v$, thus $$E_+^2(v_m) \to L:= \inf\left\{E_+^2(w): w\in\Lambda_v\right\}$$ as $m\to\infty$. Fix distinct points $p_1,p_2,p_3\in S^1$ and distinct points $\bar{p}_1, \bar{p}_2,\bar{p}_3\in \Gamma$. After possibly precomposing each $v_m$ with a conformal diffeomorphism of $D$ we may assume that every $v_m$ satisfies the $3$-point condition $\trace(v_m)(p_i) = \bar{p}_i$ for $i=1,2,3$.  By \cite[Proposition 7.4]{LW15-Plateau}, the family $\{\trace(v_m): m\in\N\}$ is equi-continuous. Thus, after possibly passing to a subsequence we may assume, by the Arzel\`a-Ascoli theorem, that $\trace(v_m)$ converges uniformly to a weakly monotone parametrization $\gamma$ of $\Gamma$. Fix $x_0\in\Gamma$. By \cite[Lemma 4.11]{LW15-Plateau} we have $$\sup_{m\in\N}\left[ \int_Dd(x_0, v_m(z))^2\,dz\right]<\infty$$ and hence $(v_m)$ is a bounded sequence.

Let $X_\omega$ be an ultra-completion of $X$ such that $X$ admits a $1$-Lipschitz retraction $P\colon X_\omega\to X$. We now apply Theorem~\ref{thm:ultra-limit-Sobolev-maps}. Thus, after possibly passing to a further subsequence, the map $w(z):= [(v_m(z))]$ belongs to $W^{1,2}(D, X_\omega)$ and satisfies $\trace(w) = \iota\circ\gamma$ as well as $$E_+^2(w)\leq \lim_{m\to\infty} E_+^2(v_m) = L$$ and $$\Area(w)\leq \liminf_{m\to\infty} \Area(v_m)\leq \Area(v).$$ Since $P$ is $1$-Lipschitz the map $u:= P\circ w$ belongs to $W^{1,2}(D, X)$ and satisfies $$\trace(u) = P\circ\iota\circ\gamma = \gamma$$ as well as $E_+^2(u)\leq L$ and $\Area(u)\leq \Area(v)$. It follows that $u\in \Lambda_v$ and consequently $E_+^2(u) = L$. 
Finally, since for every biLipschitz homeomorphism $\varrho\colon\overline{D} \to\overline{D}$ we have $u\circ\varrho\in \Lambda_v$ and thus $E_+^2(u) = L \leq E_+^2(u\circ\varrho)$ we see from \cite[Theorem 6.1]{LW15-Plateau} that $u$ is $\sqrt{2}$-quasiconformal.
\end{proof}

\begin{proof}[Proof of Theorem~\ref{thm:Plateau-1-complemented}]
 Let $(v_m)\subset \Lambda(\Gamma, X)$ be an area minimizing sequence, that is, $\Area(v_m)\to L$ as $m\to\infty$, where $$L:= \inf\{\Area(v): v\in\Lambda(\Gamma, X)\}.$$ By Lemma~\ref{lem:keep-area-lower-energy}, there exists for each $m$ some $\sqrt{2}$-quasiconformal map $u_m\in\Lambda(\Gamma, X)$ with $\Area(u_m)\leq \Area(v_m)$. In particular, $(u_m)$ is still an area minimizing sequence in $\Lambda(\Gamma, X)$ and moreover satisfies $E_+^2(u_m)\leq 2\cdot\Area(u_m)$ because $u_m$ is $\sqrt{2}$-quasiconformal. Thus, the sequence $(u_m)$ has uniformly bounded energy. Arguing exactly as in the proof of Lemma~\ref{lem:keep-area-lower-energy} we obtain the existence of $u\in \Lambda(\Gamma, X)$ such that $$\Area(u) \leq \lim_{m\to\infty} \Area(u_m) = L$$ and hence $u$ is an area minimizer in $\Lambda(\Gamma, X)$. By Lemma~\ref{lem:keep-area-lower-energy}, we may assume that $u$ is moreover $\sqrt{2}$-quasiconformal. This completes the proof of the theorem.
\end{proof}

Exactly as the corresponding result in \cite{LW15-Plateau}, Theorem~\ref{thm:Plateau-1-complemented} above holds with the parametrized Hausdorff area $\Area(u)$ replaced by the parametrized area induced by a definition of volume inducing quasi-convex $2$-volume densities.

The following result shows that sometimes the Plateau problem can be solved even if a space is not $1$-comple\-mented in an ultra-completion. Compare with \cite[Remark 4.5]{Wen05}.

\bp\label{prop:c_0-Plateau-not-1-complemented}
 The Banach space $c_0$ of sequences of real numbers converging to zero, equipped with the sup-norm, is not $1$-complemented in an ultra-completion of itself. Nevertheless, the Plateau problem can be solved in $c_0$ for every Jordan curve $\Gamma\subset c_0$ which satisfies $\Lambda(\Gamma, c_0)\not= \emptyset$.
\ep

\begin{proof}
Let $\omega$ be a non-principal ultrafilter on $\N$. The Banach space $\ell^\infty$ of bounded sequences of real numbers, equipped with the sup-norm, isometrically embeds into the ultra-completion $(c_0)_\omega$ of $c_0$ via the map $x\mapsto [(\varphi_m(x))]$, where the maps $\varphi_m\colon \ell^\infty\to c_0$ are given by $$\varphi_m(x):= (x_1,\dots, x_m, 0, 0,\dots)$$ for $x = (x_1, x_2, \dots)\in\ell^\infty$. Since there is no $1$-Lipschitz retraction of $\ell^\infty$ onto $c_0$, see for example \cite[Example 1.5]{BLGeomNonFunc}, it follows that $c_0$ is not $1$-complemented in $(c_0)_\omega$. 
 
Let now $\Gamma\subset c_0$ be a Jordan curve such that $\Lambda(\Gamma, c_0)\not= \emptyset$. Since $\Gamma$ is compact there exists a $1$-Lipschitz retraction $P\colon c_0\to K$ onto some compact set $K\subset c_0$ containing $\Gamma$. Notice that $\Lambda(\Gamma, K)\not=\emptyset$ and hence, by \cite[Theorem 1.1]{LW15-Plateau}, there exists $u\in \Lambda(\Gamma, K)$ which is $\sqrt{2}$-quasiconformal and satisfies $$\Area(u) = \inf\left\{\Area(v): v\in\Lambda(\Gamma, K)\right\}.$$
Since $P$ is $1$-Lipschitz it follows that $u$, as an element of $\Lambda(\Gamma, c_0)$, is also an area minimizer in $c_0$, thus $$\Area(u) =  \inf\left\{\Area(v): v\in\Lambda(\Gamma, c_0)\right\}.$$ This concludes the proof.
\end{proof}

We leave the verification of details in the next example to the interested reader.

\begin{example}\label{ex:complete-metric-no-Plateau}
  Let $X$ be the complete metric space obtained by gluing the metric spaces $$X_n:= \left(S^1\times [0,1/n]\right) \cup \left(S_+^2\times\{1/n\}\right)$$ for $n =1, 2, \dots,$ along their common boundary $S:= S^1\times \{0\}$, where $S_+^2$ denotes the closed upper hemisphere of the standard unit sphere in $\R^3$.  Here, every $X_n$ is equipped with the natural length metric. Since the $1$-neighborhood of $S$ in $X$ retracts Lipschitzly onto $S$ it is not difficult to see that $X$ admits a quadratic isoperimetric inequality in the sense of \cite{LW15-Plateau}. We claim that the Plateau problem in $X$ cannot be solved for the Jordan curve $\Gamma:= S$ even though $\Lambda(\Gamma, X)\not=\emptyset$. We argue by contradiction and suppose there exists $u\in \Lambda(\Gamma, X)$ with minimal area and which is $\sqrt{2}$-quasiconformal. We clearly have $$\Area(u)\leq \inf_{n\in\N} \hm^2(X_n) = \hm^2(S_+^2).$$ Since $X$ admits a quadratic isoperimetric inequality it follows from \cite{LW15-Plateau} that $u$ is continuous on $\overline{D}$ and has Lusin's property (N). For topological reasons the compact image $u(\overline{D})$ has to contain one of the subsets $X_n$ and hence, by the area formula, $$\hm^2(S_+^2)< \hm^2(X_n) \leq \Area(u) \leq \hm^2(S_+^2),$$ which is impossible. This proves the claim.
 \end{example}

The arguments used in the proof of Theorem~\ref{thm:Plateau-1-complemented} yield the following result which generalizes \cite[Theorem 7.6]{LW15-Plateau} from the setting of proper metric spaces to spaces which are $1$-complemented in some ultra-completion.

\bt\label{thm:energy-min-1-complemented}
  Let $X$ be a complete metric space and $\Gamma$ a Jordan curve in $X$ such that $\Lambda(\Gamma, X)\not= \emptyset$. If $X$ is $1$-complemented in some ultra-completion of $X$ then there exists $u\in \Lambda(\Gamma, X)$ such that $$E_+^2(u) = \inf\left\{E_+^2(v): v\in\Lambda(\Gamma, X)\right\}.$$ Every such $u$ is $\sqrt{2}$-quasiconformal.
\et

The theorem furthermore holds with $E_+^2$ replaced by the Korevaar-Schoen Dirichlet energy $E^2$ from \cite{KS93} and the constant $\sqrt{2}$ replaced by $2\sqrt{2}+\sqrt{6}$. This follows as above but uses \cite[Theorem 6.8]{LW15-Plateau} instead of \cite[Theorem 6.1]{LW15-Plateau}.

\section{Dirichlet's problem}\label{sec:Dirichlet}

In this section we provide the proof of Theorem~\ref{thm:Dirichlet}, which is similar to that of Theorem~\ref{thm:Plateau-1-complemented}. 

Let $(X, d)$ be a complete metric space and $\Omega\subset \R^n$ a bounded Lipschitz domain with $n\geq 2$. Let $w\in W^{1,p}(\Omega, X)$ for some $p>1$ and let $(u_m)\subset W^{1,p}(\Omega, X)$ be an $E_+^p$-energy minimizing sequence subject to the condition $\trace(u_m) = \trace(w)$.

\bl\label{lem:bdd-from-equal-trace}
 We have $$\sup_{m\in\N}\left[\int_\Omega d(x_0, u_m(z))^p\,dz\right] <\infty$$ for some and thus every $x_0\in X$. 
\el

\begin{proof}
By \cite[Theorem 1.12.2 and Corollary 1.6.3]{KS93} the function $$h_m(z):= d(w(z), u_m(z))$$ is in $W_0^{1,p}(\Omega)$ and satisfies $\sup_m\|\nabla h_m\|_{L^p(\Omega)}<\infty$. Thus, by the Poincar\'e inequality, the sequence $(h_m)$ is bounded in $L^p(\Omega)$. Since $w\in L^p(\Omega, X)$ the lemma follows.
\end{proof}

\begin{proof}[Proof of Theorem~\ref{thm:Dirichlet}]
Let $X$, $\Omega$, $w$ and $(u_m)$ be as above. By Lemma~\ref{lem:bdd-from-equal-trace} the sequence $(u_m)$ satisfies $$\sup_{m\in\N}\left[\int_\Omega d(x_0, u_m(z))^p\,dz + E_+^p(u_m)\right] <\infty$$ for some and thus every $x_0\in X$ and hence $(u_m)$ is a bounded sequence.
Let $X_\omega$ be an ultra-completion of $X$ such that $X$ admits a $1$-Lipschitz retraction $P\colon X_\omega\to X$. After possibly passing to a subsequence, we may assume by Theorem~\ref{thm:ultra-limit-Sobolev-maps} that the map $v(z):= [(u_m(z))]$ belongs to $W^{1,p}(\Omega, X_\omega)$ and satisfies $\trace(v) = \iota\circ \trace(w)$ as well as $$E_+^p(v)\leq \lim_{m\to\infty} E_+^p(u_m).$$ Since $P$ is $1$-Lipschitz the map $u:= P\circ v$ belongs to $W^{1,p}(\Omega, X)$ and satisfies $\trace(u) = \trace(w)$ and $E_+^p(u)\leq \lim_{m\to\infty} E_+^p(u_m)$. This completes the proof.
\end{proof}

As already mentioned in the introduction, Theorem~\ref{thm:Dirichlet} holds with the Reshetnyak energy replaced by the Korevaar-Schoen Dirichlet energy $E^p$ introduced in \cite{KS93}. For this, notice that $E_+^p$ and $E^p$ are comparable and $E^p$ is also lower semi-continuous, see \cite{KS93} and \cite{LW15-Plateau}.

\def\cprime{$'$} \def\cprime{$'$} \def\cprime{$'$}


\begin{thebibliography}{10}

\bibitem{AlvT04}
J.~C. Alvarez~Paiva and A.~C. Thompson.
\newblock Volumes on normed and {F}insler spaces.
\newblock In {\em A sampler of {R}iemann-{F}insler geometry}, volume~50 of {\em
  Math. Sci. Res. Inst. Publ.}, pages 1--48. Cambridge Univ. Press, Cambridge,
  2004.

\bibitem{Amb90}
Luigi Ambrosio.
\newblock Metric space valued functions of bounded variation.
\newblock {\em Ann. Scuola Norm. Sup. Pisa Cl. Sci. (4)}, 17(3):439--478, 1990.

\bibitem{AT04}
Luigi Ambrosio and Paolo Tilli.
\newblock {\em Topics on analysis in metric spaces}, volume~25 of {\em Oxford
  Lecture Series in Mathematics and its Applications}.
\newblock Oxford University Press, Oxford, 2004.

\bibitem{BLGeomNonFunc}
Yoav Benyamini and Joram Lindenstrauss.
\newblock {\em Geometric nonlinear functional analysis. {V}ol. 1}, volume~48 of
  {\em American Mathematical Society Colloquium Publications}.
\newblock American Mathematical Society, Providence, RI, 2000.

\bibitem{BrH99}
Martin~R. Bridson and Andr{\'e} Haefliger.
\newblock {\em Metric spaces of non-positive curvature}, volume 319 of {\em
  Grundlehren der Mathematischen Wissenschaften [Fundamental Principles of
  Mathematical Sciences]}.
\newblock Springer-Verlag, Berlin, 1999.

\bibitem{BI12}
Dmitri Burago and Sergei Ivanov.
\newblock Minimality of planes in normed spaces.
\newblock {\em Geom. Funct. Anal.}, 22(3):627--638, 2012.

\bibitem{Dierkes-et-al10}
Ulrich Dierkes, Stefan Hildebrandt, and Friedrich Sauvigny.
\newblock {\em Minimal surfaces}, volume 339 of {\em Grundlehren der
  Mathematischen Wissenschaften [Fundamental Principles of Mathematical
  Sciences]}.
\newblock Springer, Heidelberg, second edition, 2010.
\newblock With assistance and contributions by A. K{\"u}ster and R. Jakob.

\bibitem{Dou31}
Jesse Douglas.
\newblock Solution of the problem of {P}lateau.
\newblock {\em Trans. Amer. Math. Soc.}, 33(1):263--321, 1931.

\bibitem{EG92}
Lawrence~C. Evans and Ronald~F. Gariepy.
\newblock {\em Measure theory and fine properties of functions}.
\newblock Studies in Advanced Mathematics. CRC Press, Boca Raton, FL, 1992.

\bibitem{HWW93}
P.~Harmand, D.~Werner, and W.~Werner.
\newblock {\em {$M$}-ideals in {B}anach spaces and {B}anach algebras}, volume
  1547 of {\em Lecture Notes in Mathematics}.
\newblock Springer-Verlag, Berlin, 1993.

\bibitem{HKST15}
Juha Heinonen, Pekka Koskela, Nageswari Shanmugalingam, and Jeremy Tyson.
\newblock {\em Sobolev spaces on metric measure spaces}, volume~27 of {\em New
  Mathematical Monographs}.
\newblock Cambridge University Press, Cambridge, 2015.

\bibitem{Kar07}
M.~B. Karmanova.
\newblock Area and co-area formulas for mappings of the {S}obolev classes with
  values in a metric space.
\newblock {\em Sibirsk. Mat. Zh.}, 48(4):778--788, 2007.

\bibitem{Kir94}
Bernd Kirchheim.
\newblock Rectifiable metric spaces: local structure and regularity of the
  {H}ausdorff measure.
\newblock {\em Proc. Amer. Math. Soc.}, 121(1):113--123, 1994.

\bibitem{KS93}
Nicholas~J. Korevaar and Richard~M. Schoen.
\newblock Sobolev spaces and harmonic maps for metric space targets.
\newblock {\em Comm. Anal. Geom.}, 1(3-4):561--659, 1993.

\bibitem{Lan13}
Urs Lang.
\newblock Injective hulls of certain discrete metric spaces and groups.
\newblock {\em J. Topol. Anal.}, 5(3):297--331, 2013.



\bibitem{LW16-harmonic}
Alexander Lytchak and Stefan Wenger.
\newblock Regularity of harmonic discs in spaces with quadratic isoperimetric
  inequality.
\newblock {\em Calc. Var. Partial Differential Equations}, 55(4):55:98, 2016.

\bibitem{LW15-Plateau}
Alexander Lytchak and Stefan Wenger.
\newblock Area {M}inimizing {D}iscs in {M}etric {S}paces.
\newblock {\em Arch. Ration. Mech. Anal.}, 223(3):1123--1182, 2017.

\bibitem{LW-param}
Alexander Lytchak and Stefan Wenger.
\newblock Canonical parametrizations of metric discs.
\newblock {\em preprint arXiv:1701.06346}, 2017.

\bibitem{LW-isoperimetric}
Alexander Lytchak and Stefan Wenger.
\newblock Isoperimetric characterization of upper curvature bounds.
\newblock {\em Acta Math.} 221 (2018), no. 1, 159--202.

\bibitem{LWY16}
Alexander Lytchak, Stefan Wenger, and Robert Young.
\newblock Dehn functions and {H}\"older extensions in asymptotic cones.
\newblock {\em preprint arXiv:1608.00082}, 2016.

\bibitem{MZ10}
Chikako Mese and Patrick~R. Zulkowski.
\newblock The {P}lateau problem in {A}lexandrov spaces.
\newblock {\em J. Differential Geom.}, 85(2):315--356, 2010.

\bibitem{Mor48}
Charles~B. Morrey, Jr.
\newblock The problem of {P}lateau on a {R}iemannian manifold.
\newblock {\em Ann. of Math. (2)}, 49:807--851, 1948.

\bibitem{Nik79}
I.~G. Nikolaev.
\newblock Solution of the {P}lateau problem in spaces of curvature at most
  {$K$}.
\newblock {\em Sibirsk. Mat. Zh.}, 20(2):345--353, 459, 1979.

\bibitem{OvdM14}
Patrick Overath and Heiko von~der Mosel.
\newblock Plateau's problem in {F}insler 3-space.
\newblock {\em Manuscripta Math.}, 143(3-4):273--316, 2014.

\bibitem{Rad30}
Tibor Rad{\'o}.
\newblock On {P}lateau's problem.
\newblock {\em Ann. of Math. (2)}, 31(3):457--469, 1930.

\bibitem{Res97}
Yu.~G. Reshetnyak.
\newblock Sobolev classes of functions with values in a metric space.
\newblock {\em Sibirsk. Mat. Zh.}, 38(3):657--675, iii--iv, 1997.

\bibitem{Wen05}
Stefan Wenger.
\newblock Isoperimetric inequalities of {E}uclidean type in metric spaces.
\newblock {\em Geom. Funct. Anal.}, 15(2):534--554, 2005.

\bibitem{Wen11-cpt}
Stefan Wenger.
\newblock Compactness for manifolds and integral currents with bounded diameter
  and volume.
\newblock {\em Calc. Var. Partial Differential Equations}, 40(3-4):423--448,
  2011.

\end{thebibliography}
\end{document}